\documentclass[10pt]{article}
\usepackage{amsmath,amssymb,accents}
\usepackage{hyperref}
\usepackage{tikz}
\usetikzlibrary{calc,arrows}

\newcommand{\status}{}

\newcommand{\file}{$\ti{\ }$/wisk/strongtight/stight.tex}
\renewcommand{\file}{}

\newcommand{\detail}[1]{\par\noi{\bf[Proof detail\ }{#1}
\hfill{\bf ]}\par\noi\hspace{-4pt}}
\renewcommand{\detail}[1]{}

\newcommand{\dis}{\displaystyle}
\newcommand{\txt}{\textstyle}

\newcommand{\noi}{\noindent}

\newcommand{\halmos}{\rule{1ex}{1.4ex}}
\def \qed {\nopagebreak{\hspace*{\fill}$\halmos$\medskip}}
\newcommand{\med}{\medskip}

\newtheorem{theorem}{Theorem}
\newtheorem{proposition}[theorem]{Proposition}
\newtheorem{corollary}[theorem]{Corollary}
\newtheorem{conjecture}[theorem]{Conjecture}
\newtheorem{lemma}[theorem]{Lemma}

\newtheorem{remark}[theorem]{Remark}

\newcommand{\bt}{\begin{theorem}}
\newcommand{\et}{\end{theorem}}
\newcommand{\bl}{\begin{lemma}}
\newcommand{\el}{\end{lemma}}
\newcommand{\bp}{\begin{proposition}}
\newcommand{\ep}{\end{proposition}}
\newcommand{\bcor}{\begin{corollary}}
\newcommand{\ecor}{\end{corollary}}
\newcommand{\br}{\begin{remark}\rm}
\newcommand{\er}{\end{remark}}
\newcommand{\bcon}{\begin{conjecture}}
\newcommand{\econ}{\end{conjecture}}

\newcommand{\be}{\begin{equation}}
\newcommand{\ee}{\end{equation}}
\newcommand{\bes}{\begin{equation*}}
\newcommand{\ees}{\end{equation*}}
\newcommand{\bee}{\begin{enumerate}}
\newcommand{\eee}{\end{enumerate}}
\newcommand{\bei}{\begin{itemize}}
\newcommand{\eei}{\end{itemize}}
\newcommand{\bea}{\begin{eqnarray}}
\newcommand{\eea}{\end{eqnarray}}
\newcommand{\beas}{\begin{eqnarray*}}
\newcommand{\eeas}{\end{eqnarray*}}
\newcommand{\ba}{\begin{array}}
\newcommand{\ea}{\end{array}}
\newcommand{\bc}{\be\begin{array}{r@{\,}c@{\,}l}}
\newcommand{\ec}{\end{array}\ee}

\newcommand{\al}{\alpha}

\newcommand{\de}{\delta}

\newcommand{\eps}{\varepsilon}
\newcommand{\la}{\lambda}
\newcommand{\La}{\Lambda}

\newcommand{\si}{\ensuremath{\sigma}}

\newcommand{\Fi}{{\cal F}}
\newcommand{\Gi}{{\cal G}}

\newcommand{\Li}{{\cal L}}

\newcommand{\Pc}{{\cal P}}

\newcommand{\Si}{{\cal S}}

\newcommand{\R}{{\mathbb R}}
\newcommand{\N}{{\mathbb N}}
\newcommand{\Z}{{\mathbb Z}}

\newcommand{\E}{{\mathbb E}}
\newcommand{\I}{{\mathbb I}}
\renewcommand{\P}{{\mathbb P}}

\newcommand{\sub}{\subset}

\newcommand{\asto}[1]{\underset{{#1}\to\infty}{\longrightarrow}}
\newcommand{\Asto}[1]{\underset{{#1}\to\infty}{\Longrightarrow}}

\newcommand{\dgg}{\dagger}

\newcommand{\un}{\underline}
\newcommand{\subb}[2]{_{\ba{c}\scriptstyle{#1}\\[-.15cm]\scriptstyle{#2}\ea}}

\newcommand{\ffrac}[2]{{\textstyle\frac{{#1}}{{#2}}}}

\newcommand{\nab}{\nabla}

\newcommand{\half}{{[0,\infty)}}

\newcommand{\quand}{\quad\mbox{and}\quad}

\newcommand{\ha}{\ffrac{1}{2}}
\newcommand{\haa}{\frac{1}{2}}

\makeatletter\@addtoreset{equation}{section}
\makeatother 

\newcommand{\Pssi}{\Psi}

\begin{document}

\title{\vspace{-3cm}Noninvadability implies noncoexistence for a class of
 cancellative systems}
\author{Jan M. Swart\vspace{6pt}\\
{\small \' UTIA}\\
{\small Pod vod\'arenskou v\v e\v z\' i 4}\\
{\small 18208 Praha 8}\\
{\small Czech Republic}\\
{\small e-mail: swart@utia.cas.cz}\vspace{4pt}}
\date{{\scriptsize\file}\quad\today}
\maketitle\vspace{-.7cm}
\status

\begin{abstract}\noi
There exist a number of results proving that for certain classes of
interacting particle systems in population genetics, mutual invadability of
types implies coexistence. In this paper we prove a sort of converse statement
for a class of one-dimensional cancellative systems that are used to
model balancing selection. We say that a model exhibits strong interface
tightness if started from a configuration where to the left of the origin all
sites are of one type and to the right of the origin all sites are of the
other type, the configuration as seen from the interface has an invariant law
in which the number of sites where both types meet has finite expectation. We
prove that this implies noncoexistence, i.e., all invariant laws of the
process are concentrated on the constant configurations. The proof is based on
special relations between dual and interface models that hold for a large
class of one-dimensional cancellative systems and that are proved here
for the first time.
\end{abstract}

\vspace{.4cm}
\noi
%
{\it Acknowledgments.} Work sponsored by GA\v{C}R grant: P201/10/0752.


\section{Introduction and main result}

In spatial population genetics, one often considers interacting particle systems
where each site in the lattice can be occupied by one of two different types,
respresenting different genetic types of the same species or even different
species. It is natural to conjecture that if each type is able to invade an area
that is so far occupied by the other type only, then coexistence should be
possible, i.e., there should exist invariant laws that are concentrated on
configurations in which both types are present. There exist a number of
rigorous results of this nature. In particular, Durrett \cite{Dur02} has
proved a general result of this sort for systems with fast stirring; see also,
e.g., \cite{DN97} for similar results. In a more restricted context, the same
idea (mutual invadability implies coexistence) is also behind the proofs of,
e.g., \cite[Thm~1~(b)]{NP99} or \cite[Thm~4]{CP07}.

In this paper, we will prove a converse claim. We will show that for a class
of one-dimensional cancellative systems that treat the two types
symmetrically, mutual non-invadability implies noncoexistence. In particular,
this applies to several generalizations of the standard, one-dimensional voter
model that are used to model {\em balancing selection} (sometimes also called
heterozygosity selection or negative frequency dependent selection), which is
the effect, observed in many natural populations, that types that are locally
in the minority have a selective advantage, since they are able to use
resources not available to the other type.

Since our general theorem needs quite a bit of preparation to formulate, as a
warm-up and motivation for what will follow, we first describe three
particular models that our result applies to. These models also occur in
\cite{SS08hv}. We refer to that paper for a more detailed motivation and a
proof that they are indeed cancellative systems.

Restricting ourselves to the one-dimensional case, as we will throughout the
paper, let $\{0,1\}^\Z$ be the space of configurations $x=(x(i))_{i\in\Z}$ of
zeros and ones on $\Z$. We sometimes identify sets with indicator functions
and write $|x|=|\{i:x(i)=1\}|$ for the number of ones in a configuration
$x$. We recall that an interacting particle system (with two types), in one
dimension, is a Markov process $X=(X_t)_{t\geq 0}$ with state space
$\{0,1\}^\Z$ that is defined by its local transition rates \cite{Lig85}. Let
us say that such an interacting particle system is {\em type-symmetric} if its
dynamics are symmetric under a simultaneous interchanging of all types, that
is, the transition $x\mapsto x'$ happens at the same rate as the transition
$(1-x)\mapsto(1-x')$.

The first model that our result applies to is the neutral Neuhauser-Pacala
model, which is a special case of the model introduced in \cite{NP99}.
Fix $R\geq 2$ and for each $x\in\{0,1\}^\Z$, let us write
\be
f_\tau(x,i):=\frac{1}{2R}\!\sum\subb{j\in\Z}{0<|i-j|\leq R}\!1_{\{x(j)=\tau\}}
\qquad\big(\tau=0,1,\ x\in\{0,1\}^\Z,\ i\in\Z)
\ee
for the local frequency of type $\tau$ near $i$. Then the {\em neutral
 Neuhauser-Pacala model} with {\em competition parameter} $0\leq\al\leq 1$ is
the type-symmetric interacting particle system such that
\be\label{NP}
x(i)\mbox{ flips }0\mapsto 1\mbox{ with rate }
f_1(x,i)\big(f_0(x,i)+\al f_1(x,i)\big),
\ee
and similarly for flips $1\mapsto 0$, by type-symmetry. Similarly, the
{\em affine voter model} with competition parameter $0\leq\al\leq 1$ is
the type-symmetric interacting particle system such that
\be\label{affine}
x(i)\mbox{ flips }0\mapsto 1\mbox{ with rate }
\al f_1(x,i)+(1-\al)1_{\txt\{f_1(x,i)>0\}}.
\ee
The affine voter model interpolates between the threshold voter model
(corresponding to $\al=0$) studied in, e.g., \cite{CD91,Han99,Lig94} and the
usual range $R$ voter model (for $\al=1$). The neutral Neuhauser-Pacala model
likewise reduces to a range $R$ voter model for $\al=1$. Finally, the {\em
 rebellious voter model}, introduced in\cite{SS08hv}, with competition
 parameter $0\leq\al\leq 1$ is the type-symmetric interacting particle
system such that
\be\ba{r@{}l}\label{rebel}
\dis x(i)\mbox{ flips }0\leftrightarrow 1\mbox{ with rate }
&\dis\ \ha\al
\big(1_{\txt\{x(i-1)\neq x(i)\}}+1_{\txt\{x(i)\neq x(i+1)\}}\big)\\[5pt]
&\dis+\ha(1-\al)
\big(1_{\txt\{x(i-2)\neq x(i-1)\}}+1_{\txt\{x(i+1)\neq x(i+2)\}}\big).
\ec
For $\al=1$, this model reduces to the standard nearest-neighbour voter model.

Let $X$ be a type-symmetric interacting particle system. Then, by
type-symmetry, it is easy to see that the process $Y$ defined by
\be\label{inter}
Y_t(i):=1_{\txt\{X_t(i-\ha)\neq X_t(i+\ha)\}}
\qquad(t\geq 0,\ i\in\Z+\ha)
\ee
(where $\Z+\ha:=\{i+\ha:i\in\Z\}$) is a Markov process. We call $Y$ the {\em
 interface model} of $X$. We will often say that a site $i$ is occupied by a
particle if $Y_t(i)=1$; otherwise the site is empty. Under mild assumptions on
the flip rates of $X$ (e.g.\ finite range), $Y$ is itself an interacting
particle system, where always an even number of sites flip at the same
time. Let $\un 0,\un 1\in\{0,1\}^\Z$ denote the configurations that are
constantly zero or one, respectively. If $\un 0$ (and hence by
type-symmetry also $\un 1$) is a trap for the process $X$, then, under mild
assumptions on the flip rates (e.g.\ finite range), we have that
$|Y_0|<\infty$ implies $|Y_t|<\infty$ a.s.\ for all $t\geq 0$. In particular,
this applies to all models introduced above. It is easy to see that $Y$
{\em preserves parity}, i.e., $|Y_0|\ {\rm mod}(2)=|Y_t|\ {\rm mod}(2)$
a.s.\ for all $t\geq 0$. If $|Y_0|$ is finite and odd, then we let
$l_t:=\inf\{i\in\Z+\ha:Y_t(i)=1\}$ denote the position of the left-most
particle and we let
\be\label{hatY}
\hat Y_t(i):=Y(l_t+i)\qquad(t\geq 0,\ i\in\N)
\ee
denote the process $Y$ {\em viewed from the left-most particle}. Note that
$\hat Y$ takes values in the countable state space $\hat S$ of all functions
$\hat y:\N\to\{0,1\}$ such that $|\hat y|$ is finite and odd and $\hat
y(0)=1$. Let $\de_0$ denote the unique state in $\hat S$ that contains a
single particle. We let $\hat S_{\de_0}$ denote the set of states in $\hat S$
that can be reached with positive probability from the state $\de_0$.

Following terminology first introduced in \cite{CD95}, we say that a
type-symmetric interacting particle system $X$ exhibits {\em interface
  tightness} if its corresponding interface model $\hat Y$ viewed from the
left-most particle is positive recurrent on $\hat S_{\de_0}$. In particular,
this implies that the process $\hat Y$ started from $\hat Y_0=\de_0$ spends a
positive fraction of its time in $\de_0$ and is ergodic with a unique
invariant law on $\hat S_{\de_0}$. Let $\hat Y_\infty$ be distributed
according to this invariant law. Then, by definition, we will say that $X$
exhibits {\em strong interface tightness} if $\E[|\hat Y_\infty|]<\infty$.
Strong interface tightness will be our way of rigorously formulating the idea of
`noninvadability', i.e., that neither type is able to penetrate the area
occupied by the other type.

We say that $X$ exhibits {\em coexistence} if there exists an invariant law
$\mu$ such that $\mu(\{\un 0,\un 1\})=0$, i.e., $\mu$ is concentrated on
configurations in which both types are present, and we say that $X$ {\em
 survives} if the process $X$ started with a single one (and all other sites of
type zero) satisfies $\P[X_t\neq\un 0\ \forall t\geq 0]>0$. We will prove the
following theorem.

\bt{\bf(Strong interface tightness implies noncoexistence)}\label{T:intro}
Let $X$ be either a neutral Neuhauser-Pacala model, or an affine voter model,
or a rebellious voter model, with competition parameter $0<\al\leq
1$. Assume that $X$ exhibits strong interface tightness. Then $X$ exhibits
noncoexistence.
\et

To put this into context, let us look at what is known, both rigorously and
nonrigorously, about these models. Numerical simulations for the rebellious
voter model, reported in \cite{SV10}, give the following picture. There exists
a critical parameter $\al_{\rm c}\approx 0.510\pm 0.002$ such that the process
survives and coexistence holds if and only if $\al<\al_{\rm c}$, while
interface tightness holds if and only if $\al>\al_{\rm c}$ (in particular, at
$\al=\al_{\rm c}$ one has neither survival, coexistence, nor interface
tightness). Moreover, it seems that whenever interface tightness holds, one
has strong interface tightness and in fact the probability $\P[|\hat
  Y_\infty|=(2n+1)]$ decays exponentially fast in $n$.\footnote{By contrast,
  for pure voter models of range $R\geq 2$, where strong interface tightness
  has been rigorously proved, it is known that the {\em length} of the
  interface $\sup\{i\in\N:\hat Y(i)=1\}$ has a heavy-tailed distribution with
  infinite first moment \cite[Thm~1.4]{BMSV06}.} The behaviour of the neutral
Neuhauer-Pacala model and affine voter model is supposed to be similar.

Most of these numerical `facts' are unproven but for the rebellious voter
model it has been rigorously shown that for $\al$ sufficiently close to zero
one has coexistence and no interface tightness \cite[Thm~4]{SS08hv}. It is
moreover known that coexistence is equivalent to survival
\cite[Lemma~2]{SS08hv}. It is also known rigorously that the Neuhauser-Pacala
model exhibits coexistence for $\al$ sufficiently close to zero
\cite[Thm~1~(b)]{NP99} and that the affine voter model exhibits coexistence at
$\al=0$ \cite{Lig94}. It is likely this latter result can be extended to $\al$
sufficiently small.

It is an open problem to prove either noncoexistence or interface tightness
for any of these models for any $\al<1$.\footnote{\label{f:dis}Setting $R=1$
  in either the neutral Neuhauer-Pacala model or affine voter model yields, up
  to a trivial rescaling of time, the {\em disagreement voter model} (using
  terminology from \cite{SS08hv}), which is known to exhibit noncoexistence
  and interface tightness for all $0\leq\al<1$. For the special case
  $\al=\frac{1}{2}$, moreover strong interface tightness has been proved in
  \cite[Corollary of Thm~2]{ALM92}. But, as explained in \cite{SS08hv}, this
  model has special properties that give few clues on how to prove
  noncoexistence for any of the other models.} (For $\al=1$, which corresponds
to a one-dimensional pure voter model, noncoexistence and strong interface
tightness are known.) The present result, therefore, unfortunately does not
prove anything new for these models, except that it shows that {\em if} by
some means one is able to prove strong interface tightness, then this implies
noncoexistence (and of course also interface tightness).

The rest of the paper is organized as follows. In the next section, we
formulate our general result. We introduce a class of one-dimensional
cancellative systems that will be our general framework and point out
some interesting relations between their interface models and their dual
models in the sense of cancellative systems duality. In particular, we show
that each one-dimensional, type-symmetric, cancellative system $X$
has a rather peculiar dual $X'$ that is also type-symmetric and
cancellative. This sort of duality was sort of implicit in \cite{SS08hv} but
is for the first time formally written down here. We then observe that strong
interface tightness for $X$ implies the existence of a harmonic function for
$X'$ that allows us to prove that this process dies out and hence, by duality,
that noncoexistence holds for $X$. The final section of the paper contains
proofs.

\section{Methods and further results}

\subsection{Cancellative systems}

Cancellative systems are a special class of interacting particle systems
that are linear with respect to addition modulo 2. It will be convenient to
allow the lattice to be $\I=\Z$ or $\I=\Z+\ha$. It is well-known (though for
probabilists perhaps not always at the front of their minds) that linear
spaces can be defined over any field. In particular, we may view the space
$\{0,1\}^\I$ of all functions $x:\I\to\{0,1\}$ as a linear space over the
finite field $\{0,1\}$, where the latter is equipped with addition modulo 2
(and the usual product). To distinguish this from the usual addition in $\R$
(which we will sometimes also need), we will use the symbol $\oplus$ for
(componentwise) addition modulo 2.

We equip $\{0,1\}^\I$ with the product topology and let $\Li(\I)$ denote the
space of all continuous linear maps $A:\{0,1\}^\I\to\{0,1\}^\I$. The matrix
$(A(i,j))_{i,j\in\I}$ of such a linear operator is defined as
\be\label{matdef}
A(i,j):=(A\de_j)(i)\quad\mbox{where}\quad
\de_j(i):=1_{\{i=j\}}\qquad(i,j\in\I).
\ee
It is not hard to see that the continuity of $A$ is equivalent to the
requirement that $|\{j\in\I:A(i,j)=1\}|<\infty$ for all $i\in\I$ and that
\be\label{Ax}
Ax(i)=\bigoplus_{j\in\I}A(i,j)x(j)\qquad(i\in\I),
\ee
where the infinite sum reduces to a finite sum and hence is well-defined.
Identifying sets with indicator functions as we sometimes do, we associate $A$
with the set $\{(i,j):A(i,j)=1\}\sub\I^2$. We call $\Li_{\rm
  loc}(\I):=\{A\in\Li(\I):|A|<\infty\}$ (where $|A|$ denotes the cardinality
of $A\sub\I^2$) the set of {\em local} operators on $\{0,1\}^\I$

\detail{If $A$ is a matrix such that $|\{j\in\I:A(i,j)=1\}|<\infty$ for all
  $i\in\I$, then it is easy to see that (\ref{Ax}) defines a continuous linear
  map $A:\{0,1\}^\I\to\{0,1\}^\I$. Conversely, let $A$ be such a continuous
  linear map, let $A(i,j)$ be its matrix as defined in (\ref{matdef}) and
  assume that $|\{j\in\I:A(i,j)=1\}|=\infty$ for some $i\in\I$. Then we may
  choose $j_n$, all different, such that $A(i,j_n)=1$. Let $x_n(j):=1$ if
  $j=j_1,\ldots,j_n$ and $x_n(j):=0$ for all other $j$. Then, by linearity,
  $Ax_n(i)=(-1)^n$ so $Ax_n$ does not converge to a pointwise limit as
  $n\to\infty$, while $x_n$ does, contradicting the continuity of
  $A$. Finally, if $A:\{0,1\}^\I\to\{0,1\}^\I$ is a continuous linear map, and
  its matrix as defined in (\ref{matdef}) satisfies
  $|\{j\in\I:A(i,j)=1\}|<\infty$ for all $i\in\I$, then by linearity and
  continuity, it is easy to see that $Ax$ is given by (\ref{Ax}).}

Slightly specializing\footnote{In \cite{Gri79}, a cancellative system is
  defined by a percolation substructure $\Pc(\la;V,W)$ where $\la_{i,x}$ is
  the rate of the $i$-th alarm clock at a site $x\in\Z^d$, $V_{i,x}$ is the
  set of sites where particles are spontaneously born if this alarm clock
  rings, and $W_{i,x}$ is a collection of percolation arrows. We will restrict
  ourselves to the case that $V_{i,x}=\emptyset$ for all $i,x$. The matrix
  $1\oplus A$ (where $1$ is the identity matrix) corresponds to the set of
  arrows $W_{i,x}$ in \cite{Gri79}.} from the
set-up in \cite{Gri79}, we will say that an interacting particle system on
$\I$ is {\em cancellative} if for each $A\in\Li_{\rm loc}(\I)$, there exists a
rate $r(A)\geq 0$ (possibly zero), such that $X$ makes the transition
\be\label{cancel}
x\mapsto x\oplus Ax\qquad\mbox{with rate }r(A).
\ee
We will always assume that the rates are translation invariant, i.e.,
\be\label{trans}
r(A)=r(T_k(A))\quad(k\in\Z)\quad\mbox{where}
\quad T_k(A):=\{(i+k,j+k):(i,j)\in A\},
\ee
For technical convenience, we will also assume that our models are
finite range, i.e., there exists an $R<\infty$ such that
\be\label{finran}
r(A)=0\quad\mbox{whenever}\quad\exists(i,j)\in A\mbox{ with }|i-j|>R.
\ee
It follows from standard results \cite[Thm~I.3.9]{Lig85} that any collection
of rates $(r(A))_{A\in\Li_{\rm loc}(\I)}$ satisfying (\ref{trans}) and
(\ref{finran}) corresponds to a well-defined $\{0,1\}^\I$-valued Markov
process $X$.  We refer to \cite{SS08hv} for the not immediately obvious fact
that the neutral Neuhauser-Pacala model, the affine voter model, and the
rebellious voter model are cancellative systems.

It is not hard to see that a cancellative system is type-symmetric
if and only if its rates satisfy $r(A)=0$ unless
\be\label{sf}
|\{j\in\I:(i,j)\in A\}|\mbox{ is even for all }i\in\I.
\ee
Similarly, a cancellative system is parity preserving
if and only if its rates satisfy $r(A)=0$ unless
\be\label{pp}
|\{i\in\I:(i,j)\in A\}|\mbox{ is even for all }j\in\I.
\ee
We let $\Li_{\rm ts}(\I)$ and $\Li_{\rm pp}(\I)$ denote the sets of all
$A\in\Li_{\rm loc}(\I)$ satisfying (\ref{sf}) and (\ref{pp}), respectively
(where the subscript $\rm ts$ and $\rm pp$ stand for type-symmetry and parity
preservation, respectively).

\subsection{Dual and interface models}\label{S:dual}

We set
\bc
\dis S_\pm(\I)&:=&\dis
\big\{x\in\{0,1\}^\I:\lim_{i\to\pm\infty}x(i)=0\big\},\\[5pt]
\dis S_{\rm fin}(\I)&:=&\dis\big\{x\in\{0,1\}^\I:|x|<\infty\big\}
=S_-(\I)\cap S_+(\I).
\ec
If $X$ is a cancellative system, then it is not hard to check that
\be\label{infx}\left.\ba{rlll}
{\rm(i)}&\dis\E\big[\inf X_0\big]>-\infty\quad&\mbox{implies}\quad
&\E\big[\inf X_t\big]>-\infty,\\[5pt]
{\rm(ii)}&\E\big[|X_0|\big]<\infty\quad&\mbox{implies}\quad&\E[|X_t|]<\infty,
\ea\right\}\qquad(t\geq 0),
\ee
where we notationally identify sets and indicator functions as before,
i.e., $\inf x=\inf\{i\in\I:x(i)=1\}$. It follows that $X_0\in S_-(\I)$
a.s.\ implies $X_t\in S_-(\I)$ a.s.\ for all $t\geq 0$ and by symmetry
analogue statements hold for $S_+(\I)$ and $S_{\rm fin}(\I)$.

We let $xy$ denote the pointwise product of $x,y\in\{0,1\}^\I$ and write
\be
\|x\|:=\bigoplus_{i\in\I}x(i)=|x|\ {\rm mod}(2)
\qquad\big(x\in S_{\rm fin}(\I)\big).
\ee
Let $\Gi(\I,\I)$ be the set of all pairs $(x,y)$ satisfying any of the
following conditions: 1.\ $x\in S_-(\I)$ and $y\in S_+(\I)$, or 2.\ $x\in
S_+(\I)$ and $y\in S_-(\I)$, or 3.\ $x\in S_{\rm fin}(\I)$, or 4.\ $y\in
S_{\rm fin}(\I)$. We observe that the bilinear form
\be\label{inpro}
\Gi(\I,\I)\ni(x,y)\mapsto\|xy\|
\ee
is very much like an inner product. In particular,
$\|xy\|=0$ for all $y\in S_{\rm fin}(\I)$ implies $x=0$.

Let $A^\dgg(i,j):=A(j,i)$ denote the adjoint of a matrix $A$.  It follows from
general theory (see \cite[Thm~III.1.5]{Gri79}) that the cancellative system
$X$ defined by rates $(r_X(A))_{A\in\Li_{\rm loc}(\I)}$ is dual to the
cancellative system $Y'$ defined by the rates
\be\label{rdual}
r_{Y'}(A):=r_X(A^\dgg)\qquad\big(A\in\Li_{\rm loc}(\I)\big),
\ee
in the sense that
\be\label{dual}
\E\big[\|X_0Y'_t\|\big]=\E\big[\|X_tY'_0\|\big]\qquad(t\geq 0)
\ee
whenever $X$ and $Y'$ are independent (with arbitrary initial laws) and
$(X_0,Y_0)\in\Gi(\I,\I)$ a.s. Note that $\E[\|X_0Y'_t\|]=\P[|X_tY'_0|\mbox{ is
    odd}]$. By (\ref{sf}) and (\ref{pp}), $Y'$ is parity preserving if and
only if $X$ is type-symmetric. The duality in (\ref{dual}) is the analogue of
linear systems duality (see \cite[Thm~IX.1.25]{Lig85}) with normal addition
replaced by addition modulo 2.

We next consider interface models. Let us define an `interface operator' or
`discrete differential operator' $\nab:\{0,1\}^\I\to\{0,1\}^{\I+\haa}$ by
\be\label{psidef}
(\nab x)(i)=x(i-\ha)\oplus x(i+\ha)\qquad(i\in\I+\ha).
\ee
Note that if $X=(X_t)_{t\geq 0}$ is a type-symmetric cancellative 
system on $\I$ then $Y:=(\nab(X_t))_{t\geq 0}$ is its interface model as in
(\ref{inter}). Recall the definitions of $\Li_{\rm ts}(\I)$ and $\Li_{\rm
  pp}(\I)$ from (\ref{sf}) and (\ref{pp}). The next lemma says that the
interface model of each type-symmetric cancellative system is a
parity preserving cancellative system, and conversely, each parity
preserving cancellative system is the interface model of a unique
type-symmetric cancellative system.

\bl{\bf(Interface model)}\label{L:intface}
There exists a unique bijection $\Pssi:\Li_{\rm ts}(\I)\to\Li_{\rm pp}(\I+\ha)$
such that
\be\label{Psidef}
\nab Ax=\Pssi(A)\nab x\qquad\big(x\in\{0,1\}^\I\big).
\ee
Moreover, if $X$ is a type-symmetric cancellative system on $\I$
defined by rates $r_X(A)$ with $A\in\Li_{\rm ts}(\I)$, then its interface model
is the parity preserving cancellative system on $\I+\ha$ with rates
defined by
\be\label{rint}
r_Y(A):=r_X\big(\Pssi^{-1}(A)\big)\qquad\big(A\in\Li_{\rm pp}(\I+\ha)\big).
\ee
\el
An explicit formula for $\Pssi(A)$ is given in (\ref{duint}) below.

We have just seen that every type-symmetric cancellative system $X$
gives in a natural way rise to two (in most cases different) parity preserving
cancellative systems: its dual $Y'$ in the sense of (\ref{dual}) and its
interface model $Y$ as in (\ref{inter}). Now, by Lemma~\ref{L:intface}, $Y'$ is
itself the interface model of some type-symmetric cancellative 
system $X'$ and $Y$ is the dual of some type-symmetric cancellative 
system $X''$, so it seems as if continuing in this way, one could in principle
generate infinitely many different models. It turns out that this is not the
case, however. As the next lemma shows, we have $X'=X''$ and the process stops
here.

\bl{\bf(Duals and interface models)}\label{L:dualint}
Let $\Pssi:\Li_{\rm ts}(\I)\to\Li_{\rm pp}(\I+\ha)$ be as in
Lemma~\ref{L:intface}. Then
\be
\Pssi(A)^\dgg=\Pssi^{-1}(A^\dgg)\qquad\big(A\in\Li_{\rm ts}(\I)\big).
\ee
\el

Lemma~\ref{L:dualint}, together with formulas (\ref{rdual}) and (\ref{rint}),
shows that for any type-symmetric cancellative system $X$, there
exists another type-symmetric cancellative system $X'$ as well as
parity preserving cancellative systems $Y$ and $Y'$ such that the
following commutative diagram holds:
\be\label{diagram}
\begin{tikzpicture}[auto,scale=1]
\node (Yc) at (0,0) {$Y'$};
\node (Xc) at (2,0) {$X'$};
\node (X) at (0,1.3) {$X$};
\node (Y) at (2,1.3) {$Y$};
\draw[->] (X) to node{interface} (Y);
\draw[->] (Xc) to node{interface} (Yc);
\draw[<->] (X) to node[swap]{dual} (Yc);
\draw[<->] (Y) to node{dual} (Xc);
\end{tikzpicture}
\ee
An example of such a commutative diagram was given in \cite{SS08hv}, but as
far as we know, the general case is proved for the first time here. If $X$
and $X'$ are as in (\ref{diagram}), then $X$ and $X'$ are in fact themselves
dual in the following sense.

\bl{\bf(Duality of type-symmetric cancellative systems)}\label{L:XXdual}
Let $X$ be a type-symmetric cancellative system as in
(\ref{cancel})--(\ref{finran}) and let $X'$ be the dual of
the interface model of $X$, or equivalently (by Lemma~\ref{L:dualint}), let
$X$ be the dual of the interface model of $X'$. Then $X$ and $X'$ are dual in
the sense that
\be\label{Hdual}
\E\big[H(X_t,X'_0)\big]=\E\big[H(X_0,X'_t)\big]\qquad(t\geq 0)
\ee
whenever $X$ and $X'$ are independent and satisfy
$(X_0,X'_0)\in\Gi(\I,\I+\ha)$ a.s. (with $\Gi(\I,\I+\ha)$ defined analogously
to $\Gi(\I,\I)$), and $H(x,x')$ is the duality function
\be\label{Hdef}
H(x,x'):=\|(\nab x)x'\|=\|x(\nab x')\|.
\ee
\el

Using the graphical representation of cancellative systems \cite{Gri79},
the duality in (\ref{Hdual}) can be made into a strong pathwise duality. (For
this concept, and more general theory of Markov process duality, see
\cite{JK12}.)\med

\noi
{\bf Remark~1} A special property of the rebellious voter model, that in fact
motivated its introduction in \cite{SS08hv}, is that it is self-dual with
respect to the duality in (\ref{Hdual}).\med

\noi
{\bf Remark~2} It is possible for a cancellative system to be both
type-symmetric and parity preserving. In particular, this applies to the
symmetric exclusion process $Y$, which is part of a commutative diagram of the
form:\footnote{Here $X$ has pure `disagreement' dynamics (see
  footnote~\ref{f:dis}). Its dual $Z$ is an annihilating particle model, known
  as the double branching annihilating process, where particles with a certain
  rate give birth to two new particles, situated on their neighbouring
  positions.}
\be\label{exclus}
\begin{tikzpicture}[auto,scale=1]
\node (Z) at (0,0) {$Z$};
\node (Y) at (2,0) {$Y$};
\node (X) at (4,0) {$X$};
\node (Xu) at (0,1.3) {$X$};
\node (Yu) at (2,1.3) {$Y$};
\node (Zu) at (4,1.3) {$Z$};
\draw[->] (X) to node{interface} (Y);
\draw[->] (Y) to node{interface} (Z);
\draw[->] (Xu) to node{interface} (Yu);
\draw[->] (Yu) to node{interface} (Zu);
\draw[<->] (X) to node{dual} (Zu);
\draw[<->] (Y) to node{dual} (Yu);
\draw[<->] (Z) to node{dual} (Xu);
\end{tikzpicture}
\ee

\noi
{\bf Remark~3} It is interesting to speculate how much of the above goes
through if $\{0,1\}$ is replaced by a more general finite field. It seems that
at least the duality formula (\ref{dual}) holds more generally.\med

\noi
{\bf Remark~4} If $X,X',Y$ and $Y'$ are as in (\ref{diagram}), then also $Y$
and $Y'$ are dual to each other, in a sense. If $X$ and $X'$ are voter models
and $Y$ and $Y'$ are systems of annihilating random walks, then this is a form
of non-crossing duality similar to the duality of the Brownian web.

\subsection{A harmonic function}

If $X$ and $X'$ are type-symmetric cancellative systems that are
dual in the sense of (\ref{Hdual}), then it is not hard to show that
coexistence of $X$ is equivalent to survival of $X'$. In fact, this is just
\cite[Lemma~1(a)]{SS08hv}, translated into our present notation (compare also
formula (\ref{last}) below). Our strategy for proving Theorem~\ref{T:intro}
will be to show that strong interface tightness for $X$ implies extinction
of~$X'$.

It is well-known that a duality between two Markov processes translates
invariant measures of one process into harmonic functions of the other process.
Mimicking a trick used in \cite{SS11}, we will apply this to the infinite,
translation-invariant measure
\be\label{mudef}
\mu:=\sum_{i\in\I+\ha}\P\big[(\hat Y_\infty+i)\in\cdot\,\big],
\ee
where $\hat Y_\infty$ is distributed according to the invariant law of the
interface model of $X$ viewed from the left-most particle and $\hat
Y_\infty+i$ denotes the configuration obtained from $\hat Y_\infty$ by
shifting all particles by $i$. (In set-notation, $\hat
Y_\infty+i=\{j+i:j\in\hat Y_\infty\}$.) It is not hard to see that $\mu$ is
indeed an invariant measure of the interface model $Y$ of $X$. We will not
directly use this fact, but it provides the idea for the following lemma.

\bl{\bf(Harmonic function)}\label{L:harm}
Let $X$ and $X'$ be type-symmetric cancellative systems defined by rates as in
(\ref{cancel})--(\ref{finran}), on $\I$ and $\I+\ha$ respectively, that are
dual in the sense of (\ref{Hdual}). Assume that strong interface tightness
holds for $X$ and let $\hat Y_\infty$ be distributed according to the
invariant law of the interface model of $X$ viewed from the left-most
particle. Then
\be\label{hdef}
h(x):=\sum_{i\in\I+\ha}\E\big[\|(\hat Y_\infty+i)x\|\big]
\qquad\big(x\in S_{\rm fin}(\I+\ha)\big)
\ee
defines a harmonic function $h:S_{\rm fin}(\I+\ha)\to\half$ for the process
$X'$, i.e., for each deterministic initial state
$X'_0=x'\in S_{\rm fin}(\I+\ha)$, the process $M=(M_t)_{t\geq 0}$ defined by
\be\label{Mdef}
M_t:=h(X'_t)\qquad(t\geq 0)
\ee
is a martingale with respect to the filtration generated by $X'$. 
Moreover, defining constants $0<c\leq C<\infty$ by $c:=\P[|\hat Y_\infty|=1]$
and $C:=\E[|\hat Y_\infty|]$, one has that
\be\label{linbd}
c|x|\leq h(x)\leq C|x|\qquad\big(x\in S_{\rm fin}(\I+\ha)\big).
\ee
\el

We note that if $X$ is a nearest-neighbour voter model, then $X'$ is also a
nearest-neighbour voter model and $\hat Y_\infty=\de_0$ a.s. Now the harmonic
function $h$ from Lemma~\ref{L:harm} is just $h(x)=|x|$, which is a
well-known harmonic function for $X'$. Numerical simulations in \cite{SV10}
suggest that for the rebellious voter model, as $\al$ is lowered from the pure
voter case $\al=1$, the function $h$ defined in (\ref{hdef}) changes smoothly
as a function of $\al$ and can even be smoothly extended across the critical
point.

Since the process $M_t=h(X'_t)$ in (\ref{Mdef}) is a nonnegative martingale,
it converges a.s. We will show that this implies extinction for $X'$ under the
additional assumption that the dynamics of $X$ (and hence also $X'$) have a
nearest-neighbour voter component. This latter assumption is made for technical
convenience and can be relaxed; it seems however not easy to formulate simple,
sufficient, yet general conditions on the dynamics of $X$ that allow one to
conclude from the convergence of $h(X'_t)$ that $X'$ get extinct a.s. From
the a.s.\ extinction of $X'$ we obtain in fact a little more than just
noncoexistence for $X$.

\bt{\bf(Strong interface tightness implies clustering)}\label{T:clust}
Let $X$ be a type-symmetric cancellative system on $\Z$
defined by translation invariant, finite range rates as in
(\ref{cancel})--(\ref{finran}). Assume that the dynamics of $X$ have a
nearest-neighbour voter component, i.e.,
\be
r(\{(0,0),(0,1)\})\vee r(\{(0,0),(-1,0)\})>0,
\ee
and that $X$ exhibits strong interface tightness. Then, for the process started
in an arbitrary initial law,
\be\label{clust}
\P[X_t(i)=X_t(i+1)]\asto{t}1\qquad(i\in\Z).
\ee
\et

The behaviour in (\ref{clust}) is called {\em clustering} and well-known for
one-dimensional pure voter models. For pure voter models, if the initial law
of $X_0$ is translation invariant, one has moreover that
\be\label{clustp}
\P[X_t\in\cdot\,]\Asto{t}p\de_{\un 0}+(1-p)\de_{\un 1}
\quad\mbox{with}\quad
p:=\E[X_0(0)],
\ee
where $\Rightarrow$ denotes weak convergence of probability laws on
$\{0,1\}^\Z$. More generally, if $X$
satisfies the assumptions of Theorem~\ref{T:clust} and also $X'$ exhibits
interface tightness, then using duality it is not hard to check that
(\ref{clustp}) holds with
\be\label{pdef}
p:=\E\big[\|X_0\hat Y'_\infty\|\big],
\ee
where $\hat Y'_\infty$ is independent of $X_0$ and distributed according to
the invariant law of the dual of $X$ viewed from its left-most particle.  We
note that in \cite{ALM92}, clustering is proved for a variation of the range
$R$ voter model that is not a cancellative system and for which no dual is
known. The authors show that also for this model, the probability $p$ in
(\ref{clustp}) depends in a nontrivial way on the initial law, and determine
its asymptotics if $X_0$ is a product measure with low density.

\section{Proofs}

\subsection{Duality and interface models}

In this section we prove the lemmas from Section~\ref{S:dual}.

We equip $S_-(\I)$ with the stronger topology such that $x_n\to x$ if and only
if $x_n(i)\to x(i)$ for each $i\in\I$ and $\inf x_n\to\inf x$ (with notation
as in (\ref{infx})), and we let $\Li_-(\I)$ denote the space of all linear
maps $A:S_-(\I)\to S_-(\I)$ that are continuous with respect to this stronger
topology. It is not hard to see that $A\in\Li_-(\I)$ if and only if its
matrix, defined as in (\ref{matdef}), satisfies
\be\ba{l}\label{Lim}
\dis\sup\{j\in I:A(i,j)=1\mbox{ for some }i\leq k\}<\infty,\\[5pt]
\dis\inf\{i\in I:A(i,j)=1\mbox{ for some }j\geq k\}>-\infty
\ec
for all $k\in\I$. Note that for $A\in\Li_-(\I)$ and $x\in S_-(\I)$, the
infinite sum in (\ref{Ax}) reduces to a finite sum and hence is well-defined.
We define $\Li_+(\I)$ analogously. We observe that
$\Li_-(\I)\cap\Li_+(\I)\sub\Li(\I)$ and that $A\in\Li_-(\I)$ if and only if
$A^\dgg\in\Li_+(\I)$. One has
\be\label{adA}
\|x(Ay)\|=\|(A^\dgg x)y\|
\qquad\big(x\in\Si_-(\I),\ y\in\Si_+(\I),\ A\in\Li_+(\I)\big),
\ee
and the same holds if $A\in\Li_-(\I)\cap\Li_+(\I)$ and $(x,y)\in\Gi(\I,\I)$.

\detail{Proof of (\ref{Lim}): Imagine that $A$ satisfies (\ref{Lim}). Set
\bc
\dis M_k&:=&\dis\sup\{j\in I:A(i,j)=1\mbox{ for some }i\leq k\},\\[5pt]
\dis m_k&:=&\dis\inf\{i\in I:A(i,j)=1\mbox{ for some }j\geq k\}
\ec
Then $\sup\{j:A(i,j)=1\}\leq M_i<\infty$ for each $i$ and hence infinite sum
in (\ref{Ax}) reduces to a finite sum for each $x\in S_-(\I)$, so $Ax$ is
well-defined. Choose $k$ such that $x(j)=0$ for $j<k$. Then $Ax(i)=0$ for all
$i\leq m_k$ hence $Ax\in S_-(\I)$. It is clear that $A$ is linear. Assume that
$x_n\to x$ and $\inf x_n\to\inf x$. (With $\inf\un 0:=\infty$.) Set
$K:=\inf_n\inf x_n$, which satisfies $\inf x\geq K>-\infty$ by the fact that
$\inf x_n\to\inf x$. Then
\[
Ax_n(i)=\bigoplus_{j\in\I}A(i,j)x_n(j)
=\bigoplus_{K\leq j\leq M_i}A(i,j)x_n(j)
\asto{n}Ax(i)
\]
for all $i\in\I$. Moreover,
\[
\inf_n\inf Ax_n\geq m_K>-\infty
\]
which shows that $\inf Ax_n\to\inf Ax$. Conversely, if $A:S_-(\I)\to S_-(\I)$
is continuous, then it is easy to see that its matrix defined in
(\ref{matdef}) must satisfy
\[
\sup\{j\in I:A(i,j)=1\}<\infty\qquad(i\in\I),
\]
since otherwise for some $i$ we can find $j_n\to\infty$ such that
$A(i,j_n)=1$. Then $\de_{j_n}\to\un 0$ and $\inf\de_{j_n}\to\inf\un 0$ but
$A\de_{j_n}$ does not converge, contradicting continuity. It follows that
(\ref{Ax}) reduces to a finite sum for each $x\in S_-(\La)$ and equals
$Ax(i)$.

Now imagine that $m_k=-\infty$ for some $k$. Then we can find $i_n\to-\infty$
and $j_n\geq k$ such that $A(i_n,j_n)=1$. Consider the nonempty sets
$J_n:=\{j\geq k:A(i_n,j)=1\}$, which are finite by what we have already
proved. Then either there is some $j$ that is contained in infinitely many
$J_n$'s or by going to a subsequence we can assume that that $J_n$'s are
disjoint. In either case, it is easy to see that there exists some $x\in
S_-(\I)$ such that $Ax\not\in S_-(\I)$.

Finally, imagine that $M_k=\infty$ for some $k$. Since $\sup\{j\in
I:A(i,j)=1\}<\infty$ for each $i\in\I$, this is only possible if there exist
$i_n\to-\infty$ and $j_n\to\infty$ such that $A(i_n,j_n)=1$. By choosing $n$
large enough, this implies that for any $k$ there exist $i<m_k$ and $j\geq k$
such that $A(i,j)=1$, contradicting the definition of $m_k$.
}

\detail{Formula (\ref{adA}) can easily be proved by writing both expressions
 as a double sum, which reduces to a finite sum in which we can change the
 summation order.}

We define the spaces $\Li(\I,\I+\ha)$ and $\Li_\pm(\I,\I+\ha)$ of continuous
linear maps from $\{0,1\}^\I$ to $\{0,1\}^{\I+\ha}$ or from $S_\pm(\I)$ to
$S_\pm(\I+\ha)$ analogous to $\Li(\I)$ and $\Li_\pm(\I)$, respectively.
Recall the definition of the interface operator $\nab$ from (\ref{psidef}).
It is straightforward to check the following facts.

\bl{\bf(Differential operator)}\label{L:diff}
The map $\nab:S_\pm(\I)\to S_\pm(\I+\ha)$ is a bijection with inverse
$\nab^{-1}_\pm\in\Li_\pm(\I+\ha,\I)$ given by
\be\label{phidef}
\nab^{-1}_-(i,j)=1_{\{i>j\}}\quand\nab^{-1}_+(i,j)=1_{\{i<j\}}
\qquad(i\in\I,\ j\in\I+\ha).
\ee
One has $\nab^\dgg=\nab$ and $(\nab^{-1}_-)^\dgg=\nab^{-1}_+$.
\el

\noi
{\bf Remark~1} If one defines right and left discrete derivatives as
$\nab_{\rm left}x(i):=x(i+1)\oplus x(i)$ and $\nab_{\rm right}x(i):=x(i)\oplus
x(i-1)$, then $(\nab_{\rm right})^\dgg=\nab_{\rm left}$. The main reason why
we work with half-integers is that we want the operator $\nab$ to look as much
as possible like a self-adjoint operator. (Note that since $\nab$ maps
$S_\pm(\I)$ into the different space $S_\pm(\I+\ha)$, it is not strictly
speaking self-adjoint.) Half-integers are also quite natural in view of the
interpretation of $\nab$ as an interface operator.\med

\noi
{\bf Remark~2} We observe from (\ref{phidef}) and (\ref{Ax}) that
\[
\nab^{-1}_-x(i)=1_{\txt\{x(j)=1\mbox{ for an odd number of sites }j<i\}},
\]
and a similar formula holds for $\nab^{-1}_+$. Let
\be
S_{\rm even}(\I):=\dis\big\{x\in S_{\rm fin}(\I):|x|\mbox{ is even}\big\}
\quand
S_{\rm odd}(\I):=\big\{x\in S_{\rm fin}(\I):|x|\mbox{ is odd}\big\}.
\ee
Then $\nab:S_{\rm fin}(\I)\to S_{\rm even}(\I+\ha)$ is a bijection with
inverse $\nab^{-1}_-=\nab^{-1}_+$ on $S_{\rm even}(\I+\ha)$. On the other hand,
$\nab^{-1}_-x=\nab^{-1}_+x\oplus\un 1$ for $x\in S_{\rm odd}(\I+\ha)$.\med

\noi
{\bf Proof of Lemmas~\ref{L:intface} and \ref{L:dualint}} For $A\in\Li_{\rm
  pp}(\I+\ha)$, we define $\Pssi^{-1}(A)$ by
\be\label{Psinv}
\Pssi^{-1}(A):=\nab^{-1}_-A\nab=\nab^{-1}_+A\nab\qquad\big(A\in\Li_{\rm pp}(\I+\ha)\big).
\ee
Note that since $A(\,\cdot\,,j)\in S_{\rm even}(\I)$ for each $j\in\I+\ha$, in
view of Remark~2 below Lemma~\ref{L:diff}, we have
$\nab^{-1}_-A(i,j)=\nab^{-1}_+A(i,j)$ for each $i,j$, so the two formulas for
$\Pssi^{-1}(A)$ coincide. Since $\nab(i,\,\cdot\,)\in S_{\rm even}(\I)$ for
each $i\in\I+\ha$, we have that $\Pssi^{-1}(A)\in\Li_{\rm
  ts}(\I)$. Next, for $A\in\Li_{\rm ts}(\I)$, we set
\be\label{duint}
\Pssi(A):=\big(\Pssi^{-1}(A^\dgg)\big)^\dgg
=(\nab^{-1}_\mp A^\dgg\nab)^\dgg
=\nab A\nab^{-1}_\pm,
\ee
where we have used that by Lemma~\ref{L:diff}, $\nab^\dgg=\nab$ and
$(\nab^{-1}_\mp)^\dgg=\nab^{-1}_\pm$. Since $A\in\Li_{\rm ts}(\I)$ if and only
if $A^\dgg\in\Li_{\rm pp}(\I)$, this clearly defines a map $\Pssi:\Li_{\rm
  ts}(\I)\to\Li_{\rm pp}(\I+\ha)$. Now
\be
\Pssi(\Pssi^{-1}(A))x=\nab\nab^{-1}_\pm A\nab\nab^{-1}_\pm x=Ax
=\nab^{-1}_\pm\nab A\nab^{-1}_\pm\nab x=\Pssi^{-1}(\Pssi(A))x
\qquad\big(x\in S_\pm(\I+\ha)\big),
\ee
which proves that $\Pssi$ and $\Pssi^{-1}$ are each other's inverses.
Moreover,
\be
\Pssi(A)\nab x=\nab A\nab^{-1}_\pm\nab x=\nab A x\qquad\big(x\in S_\pm(\I)\big).
\ee
Since each $x\in\{0,1\}^\I$ can be written as $x=x_-\oplus x_+$ with $x_\pm\in
S_\pm(\I)$, this proves (\ref{Psidef}). Since the map
$\nab:\{0,1\}^\I\to\{0,1\}^{\I+\haa}$ is surjective, $\Pssi(A)$ is in fact
uniquely characterized by (\ref{Psidef}).

The fact that the interface model of $X$ is the parity-preserving cancellative
 system with rates as in (\ref{rint}) is immediate from (\ref{cancel}) and
(\ref{Psidef}). Lemma~\ref{L:dualint} follows from (\ref{duint}).\qed

\noi
{\bf Proof of Lemma~\ref{L:XXdual}} Immediate from Lemmas~\ref{L:intface} and
\ref{L:dualint} and cancellative systems duality (\ref{dual}). Note that the
two formulas for the duality function $H$ coincide by Lemma~\ref{L:diff}.\qed

\subsection{Noncoexistence}

{\bf Proof of Lemma~\ref{L:harm}} We start by observing that
\be\label{hup}
h(x)\leq\sum_{i\in\I+\ha}\E\big[|x(\hat Y_\infty+i)|\big]
=\sum_{i\in\I+\ha}\E\Big[\sum_{j\in x}|\de_j(\hat Y_\infty+i)|\Big]
=\E\Big[\sum_{j\in x}|\hat Y_\infty|\Big]=|x|\,\E\big[|\hat Y_\infty|\big]
\ee
for all $x\in S_{\rm fin}(\I+\ha)$, where, as we have done before, we
notationally identify $x$ with the set $\{i:x(i)=1\}$, and the second equality
is obtained by moving the sum over $i$ inside the expectation.
Similarly
\be\label{hbel}
h(x)\geq\sum_{i\in\I+\ha}
\E\big[\|x(\hat Y_\infty+i)\|1_{\{\hat Y_\infty=\de_0\}}\big]
=\P[\hat Y_\infty=\de_0]\,|x|.
\ee
Since $\E\big[|\hat Y_\infty|\big]<\infty$ and $\P[\hat Y_\infty=\de_0]>0$ by
the assumption of strong interface tightness, formulas (\ref{hup}) and
(\ref{hbel}) imply (\ref{linbd}).

The upper bound of (\ref{linbd}), together with (\ref{infx}), show that if
$\E\big[|X'_0|]<\infty$, then $\E[h(X'_t)]<\infty$ for all $t\geq 0$. Let
$Y=(Y_t)_{t\geq 0}$ be the interface model of $X$, started in the initial law
$\P[Y_0\in\cdot\,]:=\P[\hat Y_\infty\in\cdot\,]$ if $\I+\ha=\Z$ and
$\P[Y_0\in\cdot\,]:=\P[(\hat Y_\infty+\ha)\in\cdot\,]$ if $\I+\ha=\Z+\ha$, and
independent of $X'$. Then by duality (\ref{dual}), letting $l_t$ denote the
position of the left-most particle of $Y_t$, we see that
\be
\E\big[h(X'_t)\big]
=\sum_{i\in\Z}\E\big[\|X'_t(Y_0+i)\|\big]
=\sum_{i\in\Z}\E\big[\|X'_0(Y_t+i)\|\big]
=\E\big[\sum_{i\in\Z}\|X'_0(Y_t+i-l_t)\|\big]=\E\big[h(X'_0)\big],
\ee
which proves (in combination with the Markov property of $X'$) that $h(X'_t)$
is a martingale.\qed

\detail{
Private note: for the lower bound in (\ref{linbd}) we really need to use that
$\P[\hat Y_\infty=\de_0]>0$. Otherwise, one runs into the following sort of
problem. Imagine that $\hat Y_\infty=111$ a.s. and $x$ is of the form $x=\cdots
00011011011011011011011011000\cdots$. Then $\sum_i\|x(\hat Y_\infty+1)\|=2$
regardless of how large $x$ is!
}

\noi
{\bf Proof of Theorem~\ref{T:clust}} It is straightforward to check that the
one-sided nearest neighbour voter model, in which sites with rate one copy the
type on their left, is dual, in the sense of the duality in (\ref{Hdual}), to a
one-sided nearest neighbour voter model in which sites with rate one copy the
type on their right. Therefore, if the dynamics of $X$ have a (left or right)
nearest-neighbour voter component, then the the dynamics of $X'$ have a (right
or left) nearest-neighbour voter component. From this, it is easy to see that
the probability that the process $X'$ started in $x$ gets extinct
\be
q(x):=\P^x\big[\exists t\geq 0\mbox{ s.t.\ }X'_t=\un 0]
\ee
can be uniformly bounded from below in the sense that
\be\label{unext}
\inf\{q(x):|x|\leq K\}>0\qquad\forall K<\infty.
\ee
Formula (\ref{unext}) is our sole reason for assuming that the dynamics of $X$
has a (left or right) nearest-neighbour voter component; if this can be
established by some other means then the conclusions of Theorem~\ref{T:clust}
remain valid.

Extinction of $X'$ now follows from a standard argument: Letting
$(\Fi_t)_{t\geq 0}$ denote the filtration generated by $X'$, we have by the
Markov property and the a.s.\ continuity of the conditional expectation with
respect to increasing sequences of \si-fields that
\be
q(X'_t)=
\P\big[\exists s\geq 0\mbox{ s.t.\ }X'_s=\un 0\,\big|\,\Fi_t]
\asto{t}1_{\{\exists s\geq 0\mbox{ s.t.\ }X'_s=\un 0\}}\quad{\rm a.s.}
\ee
In particular, $q(X'_t)\to 0$ a.s.\ on the event that $X'$ does not get
extinct, which by (\ref{unext}) implies that $|X'_t|\to\infty$ a.s. By the
lower bound in (\ref{linbd}), it follows that $h(X'_t)\to\infty$ a.s.\ on the
event that $X'$ does not get extinct. But Lemma~\ref{L:harm} says that
$h(X'_t)$ is a nonnegative martingale, so $h(X'_t)\to\infty$ has zero
probability and hence the same must be true for the event that $X'$ does not get
extinct.

It follows that the interface model $Y'$ of $X'$ started in
$Y'_0=\de_i+\de_{i+1}$ also gets trapped in $\un 0$ a.s., so by the fact that
$Y'$ is dual to $X$ in the sense of (\ref{dual}), we find that
\be\label{last}
\P\big[X_t(i)\neq X_t(i+1)\big]
=\E\big[\|X_t(\de_i+\de_{i+1})\|\big]
=\E\big[\|X_0Y'_t\|\big]
\leq\P[Y'_t\neq\un 0]\asto{t}0.
\ee
\qed

\noi
{\bf Proof of Theorem~\ref{T:intro}} Immediate from Theorem~\ref{T:clust} and
the fact that the neutral Neuhauser-Pacala model, the affine voter model, and
the rebellious voter model are cancellative systems, which is proved in
\cite{SS08hv}.\qed

\noi
{\bf Acknowledgements} I thank two anonymous referees for bringing to my
attention the paper \cite{ALM92}, and several other helpful suggestions.


\begin{thebibliography}{Han99}

\bibitem[ALM92]{ALM92}
E.D.~Andjel, T.M.~Liggett, and T.~Mountford.
Clustering in one-dimensional threshold voter models.
{\em Stochastic Processes Appl.}~42(1), 73--90, 1992.

\bibitem[B\&06]{BMSV06}
S.~Belhaouari, T.~Mountford, R.~Sun, and G.~Valle.
Convergence results and sharp estimates for the voter model interfaces. 
{\em Electron.\ J.\ Probab.}~11, paper no.~30, 768--801, 2006.


\bibitem[CD91]{CD91}
J.T.~Cox and R.~Durrett.
Nonlinear voter models.
In {\em Random Walks, Brownian Motion and Interacting Particle Systems.
A Festschrift in Honor of Frank Spitzer}, 189--201.
Birkh\"auser, Boston, 1991.

\bibitem[CD95]{CD95}
J.T.~Cox and R.~Durrett.
Hybrid zones and voter model interfaces.
{\em Bernoulli} 1(4), 343--370, 1995.

\bibitem[CP07]{CP07}
J.T.~Cox and E.A.~Perkins.
Survival and coexistence in stochastic spatial Lotka-Volterra models.
{\em Probab.\ Theory Relat.\ Fields}~139(1-2), 89--142, 2007.

\bibitem[DN97]{DN97}
R.~Durrett and C.~Neuhauser.
Coexistence results for some competition models.
{\em Ann. Appl. Probab.}~7(1), 10--45, 1997.

\bibitem[Dur02]{Dur02}
R.~Durrett.
Mutual invadability implies coexistence in spatial models.
{\em Mem.\ Am.\ Math.\ Soc.}~740, 118 pages, 2002.

\bibitem[Gri79]{Gri79}
D.~Griffeath.
{\em Additive and Cancellative Interacting Particle Systems.}
Lecture Notes in Math.~724, Springer, Berlin, 1979.

\bibitem[Han99]{Han99}
S.J.~Handjani.
The complete convergence theorem for coexistent threshold voter models.
{\em Ann.\ Probab.}~27(1), 226--245, 1999.

\bibitem[JK12]{JK12}
S.~Jansen and N.~Kurt.
On the notion(s) of duality for Markov processes.
Preprint, 50 pages. ArXiv:1210.7193v1.

\bibitem[Lig85]{Lig85}
T.M.~Liggett.
{\em Interacting Particle Systems}.
Springer, New York, 1985.

\bibitem[Lig94]{Lig94}
T.M.~Liggett.
Coexistence in threshold voter models.
{\em Ann.\ Probab.}~22, 764--802, 1994.

\bibitem[NP99]{NP99}
C.~Neuhauser and S.W.~Pacala.
An explicitly spatial version of the Lotka-Volterra model with
interspecific competition.
{\em Ann.\ Appl.\ Probab.}~9(4), 1226--1259, 1999.

\bibitem[SS08]{SS08hv}
A.~Sturm and J.M.~Swart.
Voter models with heterozygosity selection.
{\em Ann.\ Appl. Probab.}~18(1), 59--99, 2008.


\bibitem[SS11]{SS11}
A.~Sturm and J.M.~Swart.
Subcritical contact processes seen from a typical infected site.
Preprint, 41 pages. ArXiv:1110.4777v2.

\bibitem[SV10]{SV10}
J.M.~Swart and K.~Vrbensk\'y.
Numerical analysis of the rebellious voter model.
{\em J.\ Stat.\ Phys.}~140(5), 873--899, 2010.

\end{thebibliography}
\end{document}